# CONFIDENCE INTERVALS FOR NONHOMOGENEOUS BRANCHING PROCESSES AND POLYMERASE CHAIN REACTIONS

By Didier Piau

*Université Lyon 1*

We extend in two directions our previous results about the sampling and the empirical measures of immortal branching Markov processes. Direct applications to molecular biology are rigorous estimates of the mutation rates of polymerase chain reactions from uniform samples of the population after the reaction. First, we consider nonhomogeneous processes, which are more adapted to real reactions. Second, recalling that the first moment estimator is analytically known only in the infinite population limit, we provide rigorous confidence intervals for this estimator that are valid for any finite population. Our bounds are explicit, nonasymptotic and valid for a wide class of nonhomogeneous branching Markov processes that we describe in detail. In the setting of polymerase chain reactions, our results imply that enlarging the size of the sample becomes useless for surprisingly small sizes. Establishing confidence intervals requires precise estimates of the second moment of random samples. The proof of these estimates is more involved than the proofs that allowed us, in a previous paper, to deal with the first moment. On the other hand, our method uses various, seemingly new, monotonicity properties of the harmonic moments of sums of exchangeable random variables.

**Introduction.** The incomplete replications of DNA sequences and their mutations that occur during successive cycles of a biochemical reaction called the polymerase chain reaction (PCR) can be modeled, under various simplifying hypotheses, by a branching process with a suitable branching mechanism; see Sun (1995) and Weiss and von Haeseler (1995). Sun proposed a point estimator of the mutation rate of homogeneous reactions that is valid, in fact, in the infinitely-many-sites and infinite-population limits. Sun's estimator is based on the first moment method and was adapted by Wang,

Zhang, Cheng and Sun (2000) to the finitely-many-sites case, still for the infinite-population limit of homogeneous reactions. In Piau (2002, 2004a), we showed that the branching process introduced by these authors was but an example of a wider class of processes that we called immortal branching Markov processes. We studied in-depth properties of these processes, especially in the case of polymerase chain reactions. Thus, we provided explicit bounds of the discrepancy between the point estimator of a finite-population homogeneous reaction and its infinite-population limit, in cases of both infinitely many sites and finitely many sites.

In this paper, we refine our methods and adapt them to nonhomogeneous reactions. This provides confidence intervals for the point estimator of the mutation rate. Also, we apply our results to a published data set and we comment on some estimation aspects of the model. For the sake of simplicity, we restrict the exposition to the so-called additive model, that is, to the infinitely-many-sites case, although similar results hold in the finitely-many-sites case. Finally, we show that our techniques allow us to deal with more general branching Markov processes. We explain in detail how to get pointwise estimates in this wider context and we leave as straightforward extensions the computation of confidence intervals.

Coming back to the molecular biology context, the first consequence of our results is that Sun's first moment method, supplemented by the correction that the finiteness of the initial population induces and by explicit confidence intervals, is also available for PCR with variable efficiencies. This provides an alternative to the estimation of the mutation rate through Monte Carlo simulations based on the properties of the coalescent that was proposed by Weiss and von Haeseler (1997). To our knowledge, our results are the first rigorous results that deal with nonhomogeneous reactions for finite populations. Second, we exhibit realistic efficiency sequences such that the finite-population correction is significant: In one case, we are able to show that the correct estimator is more than 33% and less than 63% higher than its infinite-population approximation for every sample. Conversely, we prove that the finite-population correction is negligible as soon as the parameters fulfill a simple condition. Third, we show that, for finite populations, the first moment method yields an estimator that is not consistent, that is, whose variance does not converge to zero when the size of the sample goes to infinity. Thus, poor confidence intervals are an intrinsic feature of this setting.

In actual reactions, the efficiency decreases along the successive cycles of the reaction (see Section 1 for a definition of the efficiency of the reaction or, more precisely, of a cycle of the reaction). The reduced sterical accessibility to the DNA sequences when the population is large is among several plausible biochemical reasons for this phenomenon. This shows that the efficiency of a cycle should be random and depend on the size of the population before that cycle. We present some extensions of our results to this setting. In particular, Schnell and Mendoza (1997) suggested that the kinetics of PCR reactions follow a Michaelis–Menten law. That is, the efficiency $\lambda_n$ of the $n$th cycle depends on the population $S_{n-1}$ before the $n$th cycle, with

$$\lambda_n = D/(C + S_{n-1}). \tag{1}$$

Here $C$ denotes the (usually quite large) Michaelis–Menten constant of the reaction, $D$ is of the order of magnitude of $C$ and $D \leq C + S_0$. (Schnell and Mendoza suggested choosing $D = C + 1$ so as to get $\lambda_1 = 1$ if $S_0 = 1$, as the greatest available efficiency.) When the initial population $S_0$ is such that $S_0 \ll C$, this allows us to recover the initial exponential growth phase, followed by a linear increase of the number of molecules; see Jagers and Klebaner (2003). Also note that Michaelis–Menten kinetics imply that when $S_0 \to \infty$, the largest value of the sequence $(\lambda_n)$, namely $\lambda_1 = D/(C + S_0)$, converges to zero. In other words, the underlying branching process becomes critical in the $S_0 \to \infty$ limit.

Point estimators and confidence intervals are consequences of precise bounds of the mean and the variance of a uniform sample. In turn, these follow from the study of the empirical measure of the population. Our methods ultimately rely on rather sharp bounds of the harmonic means of sums of i.i.d. or exchangeable random variables. Thus, on our way, we state and prove various new results about these means that are often valid in a broader context and, in particular, some simple monotonicity properties that seem to have been unnoticed until now.

The model of the PCR by a branching process is in Section 1, as well as a sample of the results of the paper. Some notation used in the paper are collected in Section 2. Theoretical results on the moments of samples are in Section 4. These follow from the results on empirical measures of Section 3. Uniform bounds are available even for random efficiencies, as explained in Section 5, and for a much more general model of branching processes, as explained in Section 6. Consequences with regard to the estimation of mutation rates are described in Section 7. In Section 8, we apply the method to the published data set used by Weiss and von Haeseler (1997). Some comments about the estimation of the efficiencies are in Section 9. Proofs are mainly deferred to Sections 10–13.

**1. Model of the PCR.** The PCR is modeled by a nondecreasing Galton–Watson process $(S_n)$ that starts from $S_0 \geq 1$ particles with a Bernoulli reproduction. We call $(S_n)$ a Bernoulli branching process. More precisely, each particle $x$ gives birth to $L_x = 1$ or $L_x = 2$ descendants independently of the other particles and with distribution

$$\mathbb{P}(L_x = 2) := \lambda =: 1 - \mathbb{P}(L_x = 1). \tag{2}$$

Each particle represents a single stranded molecule that comprises the region targeted at by the PCR or represents its complement on the other strand of the original duplex DNA molecule. Thus, the branching process counts the number $S_n$ of successfully replicated biological sequences after $n$ cycles of the reaction. Mutations are described by the states $s(x)$ of the particles $x$ as follows. Assume that the states of the $S_0$ initial particles are given. For any particle $x$ of any generation, the first descendant of $x$ is $x$ itself and it has the same state $s(x)$. If the other descendant $y$ exists, its state $s(y)$ is the result of the application of a given Markov kernel to $s(x)$. In the additive model, the states are real numbers and the kernel is given by

$$(3) \qquad s(y) := s(x) + \xi_y, \qquad \mathbb{E}(\xi_y) =: \mu, \qquad \mathbb{V}(\xi_y) =: \nu,$$

where all the random variables $\xi_y$ are independent, $\mathbb{E}$ denotes the expectation and $\mathbb{V}$ denotes the variance. In the PCR context, the law of $\xi_y$ is usually Poisson, so that $\nu = \mu$. In any case, $\lambda$ is the efficiency of the PCR (in its exponential phase) and $\mu$ is the mutation rate (per cycle and per particle). From now on, we assume that the initial population is homogeneous and we set $s(x) = 0$ for any $x$ in the initial population. Thus, in the PCR context, $s(x)$ is a nonnegative integer for any $x$ after any number of cycles. Note, however, that our results are valid in the full generality of the additive model, as described above.

In actual PCR, the efficiency $\lambda$ is not constant, but typically decreases to zero along the successive cycles of the reaction. To take into account this nonhomogeneity and the possibility of nonhomogeneous mutation rates, we choose two sequences $(\lambda_n)$ and $(\mu_n)$ indexed by $n \geq 1$, and we replace $\lambda$ and $\mu$ in (2) and (3) by $\lambda_n$ and $\mu_n$ when we construct the $n$th generation from the $S_{n-1}$ particles of the $(n-1)$th generation.

In some versions of PCR, the quantification of the product is done at the end of the reaction or only after a given number of cycles. By contrast, real-time PCR, also called quantitative PCR, allows us to measure the amount of product after each cycle. Based on fluorescent detection systems, this technology yields amplification plots that represent the accumulation of product during the successive cycles of the reaction; see Higuchi, Dollinger, Walsh and Griffith (1992) and Higuchi, Fockler, Dollinger and Watson (1993). Hence, we consider from now on that $(\lambda_n)$ is known. On the other hand, very little seems to be known about the evolution, if any, of the mutation rate during the reaction. We assume that $\mu_n = \mu$ for every $n$ and we seek to estimate the value of $\mu$. Let $\{x_1, \ldots, x_\ell\}$ denote a uniform sample of size $\ell$ drawn with replacement from the population after $n$ cycles and let $t$ denote its mean state, that is,

$$t := \ell^{-1} \sum_{i=1}^{\ell} s(x_i).$$

Because the law of $t$ is unknown, we have to rely on Bienaymé–Chebyshev bounds, which state that $t$ is in the interval bounded by

$$\mathbb{E}(t) \pm z\sqrt{\mathbb{V}(t)}$$

with probability at least $1 - 1/z^2$. This supposes known values of $\mathbb{E}(t)$ and $\mathbb{V}(t)$. Although there exists no closed form of $\mathbb{E}(t)$ and $\mathbb{V}(t)$ that would be valid for every number $n$ of generations, we can show that these quantities converge when $S_0 \to \infty$ and can compute the exact values of their limits, which we denote by $\mathbb{E}(t^*)$ and $\mathbb{V}(t^*)$, and call infinite-population limits. The task is then to bound the discrepancies between the finite-population moments and their infinite-population limits. This involves the empirical laws of the population, which are defined in Section 2 and studied in Section 3. To present a flavor of the results we are aiming at, we introduce

$$m := \mu \sum_{k=1}^{n} \frac{\lambda_k}{1 + \lambda_k}, \qquad \sigma^2 := \nu \sum_{k=1}^{n} \frac{\lambda_k}{1 + \lambda_k} + \mu^2 \sum_{k=1}^{n} \frac{\lambda_k}{(1 + \lambda_k)^2}.$$

Then, we prove that $m - \varepsilon(S_0)\mu \leq \mathbb{E}(t) \leq m$ with

$$\varepsilon(S_0) := 1/(S_0 - 1) \qquad \text{if } S_0 \geq 2, \ \varepsilon(1) := 3/2.$$

Likewise, for any $\ell \geq 3$,

$$\sigma^2/\ell \leq \mathbb{V}(t) \leq \sigma^2/\ell + (1 - 1/\ell)\eta(S_0)(\nu + \mu^2)$$

with $\eta(S_0) := 2/(S_0 - 1)$ if $S_0 \geq 2$ and $\eta(1) := 6$. Specializations of the above are the easier to establish equalities

$$\mathbb{E}(t^*) = m, \qquad \mathbb{V}(t^*) = \sigma^2/\ell.$$

Finally, we mention briefly that all these bounds on the discrepancy between the moments and the distributions of $t$ and $t^*$ are of the right order. For instance, there exists an absolute constant $c$ such that, for any $\ell \geq 3$ and $S_0 \geq 1$,

$$\mathbb{V}(t) \geq \sigma^2/\ell + c(\nu + \mu^2)/S_0.$$

An unexpected consequence is that enlarging the size of the sample becomes useless for surprisingly small sample sizes; more precisely, as soon as the deviation, which behaves like $1/S_0$, becomes the main contribution to $\mathbb{V}(t)$, instead of $\sigma^2/\ell$. Thus, $\ell \gg n^* S_0$ is useless, where $n^*$ describes the behavior of $\sigma^2$. We can choose $n^* := n$ for homogeneous reactions and $n^* := \lambda_1 + \cdots + \lambda_n$ otherwise. We recall that the expected population at time $n$ is $(1 + \lambda_1) \cdots (1 + \lambda_n) S_0$, which can be much greater than $n^* S_0$.

**2. Notation.** Call $\zeta_n$ the empirical law of the state of a particle drawn uniformly at random from the population at time $n$ and let $\eta_n := \mathbb{E}(\zeta_n)$. That is, $\zeta_n$ and $\eta_n$ are measures such that, for any nonnegative $\varphi$,

$$\zeta_n(\varphi) := S_n^{-1} \sum_{x=1}^{S_n} \varphi(s(x)), \qquad \eta_n(\varphi) := \mathbb{E}(\zeta_n(\varphi)).$$

By an abuse of notation, we denote the sum over the population at time $n$ by a sum from $x = 1$ to $x = S_n$. For any measure $\eta$, $\mathbb{M}(\eta)$ and $\mathbb{M}_2(\eta)$ denote the first and the second moments of $\eta$, and $\mathbb{D}(\eta)$ denotes its variance, that is,

$$\mathbb{M}(\eta) := \int s\, d\eta(s), \qquad \mathbb{M}_2(\eta) := \int s^2\, d\eta(s), \qquad \mathbb{D}(\eta) := \mathbb{M}_2(\eta) - \mathbb{M}(\eta)^2.$$

In the next sections, some technical lemmas are valid in the broader setting of a general nondecreasing branching process and even for the harmonic moments of sums of i.i.d. or exchangeable random variables.

DEFINITION 1. Let $L_i$ denote i.i.d. or exchangeable copies of a square integrable random variable $L$ such that $L \geq 1$ a.s. and

$$M_k := L_1 + \cdots + L_k.$$

For any $k \geq 1$, define

$$H(k) := \mathbb{E}(k/M_k), \qquad A(k) := H(k) - 1/\mathbb{E}(L), \qquad G(k) := \mathbb{E}((k/M_k)^2).$$

For any $k \geq 1$, $A(k) \geq 0$. For i.i.d. sequences, $A(k) \to 0$ when $k \to \infty$. For Bernoulli branching processes, that is, when the distribution of $L$ is given by (2), we write $A(k, \lambda)$ instead of $A(k)$ to specify the value of $\lambda$. The same convention holds for $H$ and $G$ and other sequences that are defined later.

Our results involve various parameters, functions of $(\lambda_k)$ only, that we define below.

DEFINITION 2. Let $\alpha_k := \lambda_k/(1 + \lambda_k)$, $\gamma_0 := 1$, $\gamma_0^{(i)} := 1$, and, for $n \geq 1$,

$$\gamma_n := \prod_{k=1}^{n}(1 - \alpha_k), \qquad \gamma_n^{(i)} := \prod_{k=1}^{n}(1 - \lambda_k/i).$$

Define $W_0 := 0$, $W_0' := 0$ and, for $n \geq 1$,

$$W_n := \sum_{k=1}^{n} \alpha_k, \qquad W_n' := \sum_{k=1}^{n} \alpha_k(1 - \alpha_k).$$

**3. Empirical laws.** From Theorem A, when $S_0 \to \infty$, $\eta_n$ and $\zeta_n$ converge to a deterministic measure $\eta_n^*$, which is easy to describe [we omit the proof; see Piau (2004a)].

THEOREM A. *The law $\eta_n^*$ coincides with the distribution of the random variable*

$$\varepsilon_1 \xi_1 + \cdots + \varepsilon_n \xi_n,$$

*where all the random variables $\varepsilon_k$ and $\xi_j$ are independent, $\varepsilon_k$ is Bernoulli with $\mathbb{P}(\varepsilon_k = 1) = \alpha_k = 1 - \mathbb{P}(\varepsilon_k = 0)$, and $\xi_j$ follows the law used in* (3) *at the $j$th generation. Thus,*

$$\mathbb{M}(\eta_n^*) = \sum_{k=1}^{n} \mu_k \alpha_k, \qquad \mathbb{D}(\eta_n^*) = \sum_{k=1}^{n} \nu_k \alpha_k + \mu_k^2 \alpha_k (1 - \alpha_k).$$

3.1. *Mean and distance in total variation.* The approximations below stem from precise estimations of $A(k)$ and of the harmonic moments of $S_n$ that we develop later in this paper. The results of this section are adapted from Piau (2002, 2004a) and we omit their proofs.

THEOREM B. (i) *First moment:*

$$\mathbb{M}(\eta_n) = \mathbb{M}(\eta_n^*) - \sum_{k=1}^{n} \mu_k A_k, \qquad A_k := \mathbb{E}(A(S_{k-1}, \lambda_k)).$$

(ii) *Distance in total variation:*

$$\|\eta_n - \eta_n^*\|_{\mathrm{TV}} \leq V_n := \sum_{k=1}^{n} A_k.$$

THEOREM C (Approximations). *For any $(\lambda_i)$, $S_0$ and $k \geq 1$, $A_k$ satisfies the inequalities*

$$\begin{aligned}
(S_0 + 1) A_k &\geq \gamma_{k-1} \alpha_k (1 - \lambda_k)/(1 + \lambda_k)^2, \\
(S_0 - 1) A_k &\leq \gamma_{k-1} \alpha_k (1 - \lambda_k)/(1 + \lambda_k)^2, \\
(S_0 + 1) A_k &\leq \gamma_{k-1}^{(3)} \alpha_k (1 - \lambda_k).
\end{aligned}$$

This implies that $V_n$ is bounded above and below by explicit functions of $(\mu_k)$ and $(\lambda_k)$, divided by $(S_0 - 1)$ or $(S_0 + 1)$. We assume from now on that the law of $\xi$ is constant along the generations or, more precisely, that its first two moments are. This is only for simplicity of notation and the reader should be able to guess the correct formulation of our results for variable laws of $\xi$ by analogy with the expressions in Theorem A. Thus, $\mu := \mathbb{E}(\xi)$, $\nu := \mathbb{V}(\xi)$ and

$$\mathbb{M}(\eta_n^*) = \mu W_n, \qquad \mathbb{M}(\eta_n) = \mu W_n - \mu V_n.$$

The first assertion of Theorem C provides a lower bound of $V_n$ for any $S_0$. The second assertion provides an upper bound of $V_n$ for any $S_0 \geq 2$ that involves $1/(S_0 - 1)$. When $S_0 = 1$, we should use the third assertion instead. In the rest of the paper, the bounds that involve $1/(S_0 + 1)$ are valid for any $S_0$ and the bounds that involve $1/(S_0 - 1)$ should be used only when $S_0 \geq 2$. For instance, in Corollary 3, we should use the $v_n$ lower bound for any $S_0$, the $v_n$ upper bound for any $S_0 \geq 2$ and the $v_n''$ upper bound if $S_0 = 1$ (for small values of $S_0$, the $v_n'$ bound is often less interesting than the $v_n''$ bound). These remarks apply to later results in this paper.

COROLLARY 3. (i) *One has $v_n \leq (S_0 + 1)V_n \leq v_n''$ with*

$$v_n := \sum_{k=1}^{n} \gamma_{k-1} \alpha_k (1 - \lambda_k)/(1 + \lambda_k)^2,$$

$$v_n'' := \sum_{k=1}^{n} \gamma_{k-1}^{(3)} \alpha_k (1 - \lambda_k).$$

(ii) *One has $(S_0 - 1)V_n \leq v_n$ and $S_0 V_n \leq v_n'$, with*

$$v_n' := \sum_{k=1}^{n} \gamma_{k-1}^{(2)} \alpha_k (1 - \lambda_k).$$

COROLLARY 4. *For any $S_0 \geq 2$,*

$$\mu W_n - \mu v_n/(S_0 - 1) \leq \mathbb{M}(\eta_n) \leq \mu W_n - \mu v_n/(S_0 + 1)$$

*and*

$$\|\eta_n - \eta_n^*\|_{\mathrm{TV}} \leq v_n/(S_0 - 1).$$

REMARK 5. In contrast to the last statement above, we can show that, although $\zeta_n^* = \eta_n^*$, for any given $n$ and $(\lambda_k)$, there exists a constant $c$ such that

$$\mathbb{E}[\|\zeta_n - \eta_n^*\|_{\mathrm{TV}}] \geq c/\sqrt{S_0}$$

for every $S_0 \geq 1$. We omit the proof.

3.2. *Variance and uniform bounds of the empirical laws.* We move to entirely new results, namely the estimation of second moments. Recall that we assumed for simplicity of notation that the two first moments of $\xi$ are constant. Thus,

$$\mathbb{D}(\eta_n^*) = \nu W_n - \mu^2 W_n'.$$

DEFINITION 6. Define $V_n'$ such that $0 \leq V_n' \leq V_n$ by the formula

$$V_n' := \sum_{k=1}^n A_k', \qquad A_k' := A_k^2 + (1 - 2\alpha_k) A_k.$$

The assertion $V_n' \leq V_n$ in the definition above follows from $A(\cdot, \lambda) \leq \alpha$.

THEOREM D. *One has* $\mathbb{D}(\eta_n) = \nu(W_n - V_n) + \mu^2(W_n' - V_n')$.

PROPOSITION 7. *There exists a constant $V$ that depends on $S_0$ such that*

$$\mu W_n - \mu V \leq \mathbb{E}(\eta_n) \leq \mu W_n, \qquad \|\eta_n - \eta_n^*\|_{\mathrm{TV}} \leq V,$$

$$\nu W_n + \mu^2 W_n' - (\nu + \mu^2)V \leq \mathbb{D}(\eta_n) \leq \nu W_n + \mu^2 W_n'.$$

*This holds with $V := 1/(S_0 - 1)$ for $S_0 \geq 2$, and $V := 3/2$ for $S_0 = 1$.*

3.3. *Additional term.* From Sections 3.1 and 3.2, the two first moments and the distance in total variation of the empirical laws are described by $V_n$ and $V_n'$. The second moment of uniform samples involves an additional term $R_n$, defined as

$$R_n := \mathbb{V}(\mathbb{M}(\zeta_n)) = \mathbb{E}(\mathbb{M}(\zeta_n)^2) - \mathbb{E}(\mathbb{M}(\zeta_n))^2.$$

We now complete the results of Sections 3.1 and 3.2 with an in-depth study of $R_n$. Theorem E recursively describes the evolution of $R_n$. We control the terms of the recurrence in Lemmas 9 and 11 and Corollary 14, and finally get tractable bounds of $R_n$ in Corollary 15. This section details the path that leads to these bounds of $R_n$, but the proofs of the steps themselves are postponed to Section 11.

DEFINITION 8. Introduce

$$B(k) := \mathbb{E}\left(\frac{M_k - k}{M_k^2}\right), \qquad B'(k) := \mathbb{V}\left(\frac{k}{M_k}\right), \qquad B''(k) := \frac{k}{2}\mathbb{E}\left(\frac{(L_1 - L_2)^2}{M_k^2}\right).$$

Thus, $B(k) = (H(k) - G(k))/k$ and $B'(k) = G(k) - H(k)^2$.

The theorem below is proved in Section 10.

THEOREM E. *One has $R_0 = 0$ and*

$$R_{n+1} = R_n + \nu \mathbb{E}(B(S_n)) + \mu^2 \mathbb{E}(B'(S_n)) + \mathbb{E}[\mathbb{D}(\zeta_n)B''(S_n)],$$

*where one uses $\lambda_{n+1}$ in the definition of $B(S_n)$, $B'(S_n)$ and $B''(S_n)$.*

Lemma 9 deals with the $B$ and $B'$ terms. The control of the $B''$ term is more intricate and involves Lemma 11 and Corollary 14.

LEMMA 9. *For any $k \geq 1$, $B(k) \leq b/k$ with $b := (\mathbb{E}(L) - 1)/\mathbb{E}(L)^2$. In the Bernoulli case, $b = \alpha(1-\alpha)$. For Bernoulli processes, $B'(k) \leq b'/(k+1)$ and $B''(k) \leq b''/(k+2)$, with*

$$b' := \lambda, \qquad b'' := \lambda(1-\lambda).$$

DEFINITION 10. *Let $y \geq -1$ denote a real number. Define nonnegative sequences $C$, $C'$ and $C''$ that depend on $y$ by the following equations:*

(i) *Let $C(1) := 0$ and, for any $k \geq 2$,*

$$C(k) := k^2(k+y)\mathbb{E}\left(\frac{L_1 L_2}{M_k^2(M_k+y)}\right).$$

(ii) *For any $k \geq 1$ such that $k + y > 0$,*

$$C'(k) := \mathbb{E}\left(\frac{(M_k-1)(M_k-k)}{M_k^2(M_k+y)}\right), \qquad C''(k) := \mathbb{E}\left(\frac{k(M_k-k)}{M_k^2(M_k+y)}\right).$$

In the lemma below, $\mathcal{F}_n$ is the $\sigma$-algebra generated by the $n$ first generations of the process.

LEMMA 11. *Using $\lambda_{n+1}$ in the definition of the sequences $C$, $C'$ and $C''$, we have, on the set $\{S_n + y > 0\}$,*

$$\mathbb{E}\left(\frac{\mathbb{D}(\zeta_{n+1})}{S_{n+1}+y}\Big|\mathcal{F}_n\right) = C(S_n)\frac{\mathbb{D}(\zeta_n)}{S_n+y} + C'(S_n)\nu + C''(S_n)\mu^2.$$

DEFINITION 12. *For any real number $y$ and any integer $k > -y$, let*

$$H_y(k) := \mathbb{E}\left(\frac{k+y}{M_k+y}\right).$$

Thus, $H_0(k) = H(k)$. Since $M_k \geq k$, $C''(k) \leq C'(k)$ and

$$(k+y)C'(k) \leq 1 - H(k).$$

On the other hand, for $k \geq 2$,

$$C(k) = H_y(k) - \frac{k(k+y)}{2}\mathbb{E}\bigg(\frac{(L_1 - L_2)^2}{M_k^2(M_k + y)}\bigg).$$

Thus, $C(k) \leq H_y(k)$. At this point, the only additional tools that we need are estimations of $H$ and $H_y$. These are developed in Section 12 and yield the next lemma.

LEMMA 13. (i) *For any $y \geq 0$, $C(k) \leq 1 - \lambda/(y+2)$.*
(ii) *For $y = -1$ and $k \geq 2$, $C(k) \leq 1 - \alpha$.*
(iii) *For any $y > 1 - k$, $(k+y)C'(k) \leq \alpha$.*

Our next result states that Lemma 11 can be integrated to get a recursion of the form

$$(4) \qquad \mathbb{E}\bigg(\frac{\mathbb{D}(\zeta_{n+1})}{S_{n+1} + y}\bigg) \leq \beta_{n+1}\mathbb{E}\bigg(\frac{\mathbb{D}(\zeta_n)}{S_n + y}\bigg) + \alpha_{n+1}\mathbb{E}\bigg(\frac{1}{S_n + y}\bigg)(\nu + \mu^2).$$

COROLLARY 14. *For any $y \geq 0$, (4) holds with $\beta_{n+1} := 1 - \lambda_{n+1}/(y+2)$. If $y = -1$ and $S_0 \geq 2$, (4) holds with $\beta_{n+1} := 1/(1 + \lambda_{n+1})$.*

We are now in the position to estimate the three sums that the iteration of Theorem E yields. Assume first that $S_0 \geq 2$. A weaker form of Lemma 9 is that $B(k) \leq b/(k-1)$, $B'(k) \leq b'/(k-1)$ and $B''(k) \leq b''/(k-1)$. For the $B$ and $B'$ parts, Corollary 27 gives an upper bound of $\mathbb{E}[1/(S_k - 1)]$. For the $B''$ part, we use the $y = -1$ form of Corollary 14.

Assuming now that $S_0 = 1$, we use the full form of Lemma 9 for $B$ and $B'$. Corollary 27 allows us to bound $\mathbb{E}(1/S_k)$ for the $B$ term and $\mathbb{E}[1/(S_k + 1)]$ for the $B'$ term. For the $B''$ term, we use the fact that $B''(k) \leq b''/k$ and the $y = 0$ form of Corollary 14. This yields an upper bound of $R_n$ of the form

$$R_n \leq \nu u(n) + \mu^2 u'(n) + (\nu + \mu^2)u''(n).$$

COROLLARY 15. *If $S_0 \geq 2$, we can choose*

$$u(n) := \sum_{k=1}^{n} \alpha_k(1 - \alpha_k)\gamma_{k-1}/(S_0 - 1),$$

$$u'(n) := \sum_{k=1}^{n} \lambda_k \gamma_{k-1}/(S_0 - 1),$$

$$u''(n) := \sum_{k=1}^{n-1} \lambda_k \sum_{i=k}^{n-1} \lambda_{i+1}(1 - \lambda_{i+1})\gamma_i/(S_0 - 1).$$

If $S_0 \geq 1$, we can choose

$$u(n) := \sum_{k=1}^{n} \alpha_k(1-\alpha_k)\gamma_{k-1}^{(2)}/S_0,$$

$$u'(n) := \sum_{k=1}^{n} \lambda_k \gamma_{k-1}^{(3)}/(S_0+1),$$

$$u''(n) := \sum_{k=1}^{n-1} \frac{\alpha_k}{1-\lambda_k/2} \sum_{i=k}^{n-1} \lambda_{i+1}(1-\lambda_{i+1})\gamma_i^{(2)}/S_0.$$

Uniform bounds follow from the tricks described in Section 5.

COROLLARY 16.  *There exist constants $U$, $U'$ and $U''$ that depend on $S_0$ such that $u(n) \leq U$, $u'(n) \leq U'$ and $u''(n) \leq U''$. For $S_0 \geq 2$, this holds with $U = U' = U'' = 1/(S_0 - 1)$. For $S_0 = 1$, this holds with $U = 2$, $U' = 3/2$ and $U'' = 4$.*

COROLLARY 17.  *For $S_0 \geq 2$ (resp. for $S_0 = 1$),*

$$R_n \leq 2(\nu + \mu^2)/(S_0 - 1) \qquad (\text{resp. } R_n \leq 6\nu + 11\mu^2/2).$$

**4. Moments of uniform samples.** Recall that the sample is $\{x_1, \ldots, x_\ell\}$, that the family $[s(x_i)]$ is exchangeable, and that each $s(x_i)$ follows the law $\eta_n$. Thus, first taking the expectation with respect to the randomness of the sampling procedure, and then the expectation with respect to the branching process and to the mutation process (we skip the details), we get

$$\mathbb{E}(t) = \mathbb{M}(\eta_n), \qquad \mathbb{V}(t) = \mathbb{D}(\eta_n)/\ell + (1 - 1/\ell)R_n.$$

Hence, the results below are mostly corollaries to Section 3. The exception is Proposition 19 whose proof is in Section 11 [we omit the proof of part (iv), which is anecdotal].

THEOREM F.  *One has $\mathbb{E}(t) = \mu W_n - \mu V_n$, where $V_n$ is nonnegative and converges to 0 when $S_0 \to \infty$. More precisely, for any $S_0 \geq 2$,*

$$v_n/(S_0+1) \leq V_n \leq v_n/(S_0-1)$$

*for a positive constant $v_n$, which depends on $(\lambda_k)$ only and is defined in Section 3.1.*

THEOREM G.  *We have*

$$\mathbb{V}(t) = (\nu W_n + \mu^2 W_n')/\ell - Z_n/\ell + (1 - 1/\ell)R_n,$$

where $Z_n$ and $R_n$ are nonnegative and converge to $0$ when $S_0 \to \infty$. More precisely,

$$Z_n := \nu V_n + \mu^2 V'_n \leq (\nu + \mu^2) V_n, \qquad R_n \leq r_n (\nu + \mu^2)/S_0$$

for a positive constant $r_n$ that depends on $(\lambda_k)$ only and whose value can be deduced from Corollary 15.

COROLLARY 18. *When $n \to \infty$, $\mathbb{E}(t) \to \infty$ if and only if $\mathbb{V}(t) \to \infty$ if and only if $(\lambda_k)$ is not summable. For any $(\lambda_k)$, $Z_n$ and $R_n$ are uniformly bounded.*

PROPOSITION 19.  (i) *For any $\ell \geq 1$, $\mathbb{E}(t) \leq \mathbb{E}(t^*)$.*
 (ii) *For $\ell = 1$, $\mathbb{V}(t) = \mathbb{D}(\eta_n) < \mathbb{D}(\eta_n^*) = \mathbb{V}(t^*)$.*
 (iii) *For $\ell \geq 3$, $\mathbb{V}(t) > \mathbb{V}(t^*)$.*
 (iv) *For $\ell = 2$, both situations are possible for any law of $\xi$. That is, there exist generations $n$ and efficiencies $(\lambda_k)$ such that, for any law of $\xi$, $\mathbb{V}(t) < \mathbb{V}(t^*)$, respectively, $\mathbb{V}(t) > \mathbb{V}(t^*)$.*

Finally, uniform bounds hold that are valid for any $(\lambda_k)$.

PROPOSITION 20.  (i) *For any $S_0 \geq 2$,*

$$\mu W_n - \mu/(S_0 - 1) \leq \mathbb{E}(t) \leq \mu W_n.$$

*For $S_0 = 1$, $\mu/(S_0 - 1)$ above should be replaced by $3\mu/2$.*
 (ii) *Assume that $\ell \geq 3$ and recall that $\mathbb{V}(t^*) = (\nu W_n + \mu^2 W'_n)/\ell$. For any $S_0 \geq 2$,*

$$\mathbb{V}(t^*) \leq \mathbb{V}(t) \leq \mathbb{V}(t^*) + (1 - 1/\ell) 2(\nu + \mu^2)/(S_0 - 1).$$

*For $S_0 = 1$,*

$$\mathbb{V}(t^*) \leq \mathbb{V}(t) \leq \mathbb{V}(t^*) + (1 - 1/\ell) 6(\nu + \mu^2).$$

 (iii) *Assume that $\ell = 1$. Then, for any $S_0 \geq 2$,*

$$\mathbb{V}(t^*) - (\nu + \mu^2)/(S_0 - 1) \leq \mathbb{V}(t) \leq \mathbb{V}(t^*).$$

*For $S_0 = 1$, $(\nu + \mu^2)/(S_0 - 1)$ above should be replaced by $3(\nu + \mu^2)/2$.*

**5. Random efficiencies.** Assume that $\lambda_{n+1}$ depends on $S_n$, as in the Michaelis–Menten setting we recalled in the Introduction, or on the full past $\mathcal{F}_n$ of the process up to time $n$. Perhaps surprisingly, some uniform bounds of the error term that we proved in the deterministic case still hold, but the behavior of the main term becomes somewhat unclear. In this section, we restrict the exposition to estimation of the first moment. Theorem H deals with the "error term" in the general case and Theorem I deals with the "main term" in the Michaelis–Menten case.

THEOREM H. *Let $w_n := \mathbb{E}(W_n)$, where $W_n := \alpha_1 + \cdots + \alpha_n$ is now random, and let $V$ be defined as in Proposition 7. Then*

$$\mu w_n - \mu V \leq \mathbb{E}(t) \leq \mu w_n.$$

Recall that $V \sim 1/S_0$ when $S_0 \to \infty$. Theorem H leaves open the question of the true behavior of $\mathbb{E}(t)$ in many interesting situations. For instance, the Michaelis–Menten law implies that $w_n \sim nD/S_0$ when $S_0 \to \infty$, all the other parameters being fixed. Thus, the main term $w_n$ and the error term $V$ become of the same order. Before coming back to the Michaelis–Menten case, we sketch the proof of Theorem H. We first mention without proof the crucial identities

$$\sum_{k=1}^{n} \lambda_k \gamma_k = 1 - \gamma_n, \qquad \sum_{k=1}^{n} \lambda_k \gamma_{k-1}^{(i)} = i(1 - \gamma_n^{(i)}).$$

SKETCH OF THE PROOF OF THEOREM H. A simple consequence of the monotonicity of $H$ (see Lemma 30) is

$$\mathbb{E}\left(\frac{1}{S_n} - \frac{1}{S_{n+1}} \Big| S_n\right) \geq \frac{\lambda_{n+1}}{2S_n}.$$

Taking expectations of both sides and summing over $n$, we get

$$\sum_{n \geq 0} \mathbb{E}\left(\frac{\lambda_{n+1}}{S_n}\right) \leq \frac{2}{S_0}.$$

Likewise, for any $k \geq 0$,

$$\frac{1}{S_0 - 1} - \mathbb{E}\left(\frac{1}{S_{k+1} - 1}\right) \geq \sum_{n=0}^{k} \mathbb{E}\left(\frac{\alpha_{n+1}}{S_n}\right) \geq \frac{1}{S_0} - \mathbb{E}\left(\frac{1}{S_{k+1}}\right).$$

Thus, if $(\lambda_k)$ is not summable, $S_k \to \infty$ a.s. and

$$\frac{1}{S_0 - 1} \geq \sum_{n \geq 0} \mathbb{E}\left(\frac{\alpha_{n+1}}{S_n}\right) \geq \frac{1}{S_0}.$$

These bounds are tight since $S_n = 2^n S_0$ when $\lambda_n = 1$ for every $n$. □

We now study $w_n$ in the Michaelis–Menten case, that is, when

$$w_n = \sum_{k=1}^{n} \mathbb{E}\left(\frac{\lambda_k}{1 + \lambda_k}\right) = \sum_{k=1}^{n} \mathbb{E}\left(\frac{D}{D + C + S_{k-1}}\right),$$

which depends on $n$ and $(S_0, C, D)$. Easy remarks are that $S_n \sim Dn$ almost surely and $w_n \sim \log n$ when $n \to \infty$, all the other parameters being fixed, and that $w_n \sim nD/S_0$ when $S_0 \to \infty$, all the other parameters being fixed.

Estimations for fixed values of $n$ and $S_0$ are as follows. Introduce the reduced variables

$$s_0 := S_0/C, \qquad b := C/D,$$

and note that $b(1 + s_0) \geq 1$ since $\lambda_1 = D/(C + S_0) \leq 1$. The regime we are interested in is when $s_0$ is small and $b$ is about 1, but the following result makes no such assumption.

THEOREM I. *In the Michaelis–Menten case, $w_n^- \leq w_n \leq w_n^+$ with*

$$w_n^+ := (2 + (2b - 1)/s_0) \log(1 + ns_0/(2 + s_0)),$$
$$w_n^- := \log(1 + n/(1 + b(1 + s_0))).$$

*When $b \geq 1$, we can choose $w_n^+ := w_n^*$ with*

$$w_n^* := (2 + (2b - 1)/s_0) \log(1 + ns_0/(2b(1 + s_0)^2)).$$

In the special case $b = 1$, we get

$$\log(1 + n/(2 + s_0)) \leq w_n \leq (2 + 1/s_0) \log(1 + ns_0/(2(1 + s_0)^2)).$$

PROOF OF THEOREM I. The convexity of the function $x \mapsto 1/x$ yields

$$w_{n+1} - w_n \geq D/(D + C + \mathbb{E}(S_n)).$$

Since $\mathbb{E}(S_{n+1}) = \mathbb{E}(S_n) + \mathbb{E}(S_n \lambda_{n+1})$ and $S_n \lambda_{n+1} \leq D$,

$$\mathbb{E}(S_n) \leq S_0 + nD.$$

This yields $w_n \geq \zeta(n, b(1 + s_0))$, where

$$\zeta(n, t) := \sum_{k=1}^{n} 1/(t + k) \geq \log(1 + n/(t + 1)).$$

This proves the $w_n^-$ bound. On the other hand, the concavity of the function $x \mapsto x/(1 + x)$ yields

$$w_{n+1} - w_n \leq D\mathbb{E}(1/S_n)/(1 + (D + C)\mathbb{E}(1/S_n)).$$

From Lemma 27,

$$\mathbb{E}(1/S_{n+1}|\mathcal{F}_n) \leq (1 - \lambda_{n+1}/2)/S_n.$$

For Michaelis–Menten values of $\lambda_{n+1}$, this yields

$$\mathbb{E}(1/S_{n+1}) \leq (1 - D/(2C))\mathbb{E}(1/S_n) + (D/(2C))\mathbb{E}(1/(C + S_n)).$$

The same concavity inequality with respect to $1/S_n$ that we used a few lines above allows us to deal with the term $\mathbb{E}(1/(C + S_n))$. This yields

$$\mathbb{E}(1/S_{n+1}) \leq \psi(\mathbb{E}(1/S_n))$$

for the function $\psi$ defined by

$$\psi(x) := (1 - 1/(2b))x + 1/(2b)x/(1 + Cx).$$

It is a simple matter to show that $1/\psi(x) \geq 1/a + 1/x$ for any $x \leq 1$, with

$$a := 2/D + (2b - 1)/S_0.$$

Thus $\mathbb{E}(1/S_n) \leq \psi^{(n)}(1/S_0) \leq a/(n + aS_0)$. This yields an upper bound of $w_{n+1} - w_n$ which reads, after some cumbersome algebra,

$$w_n \leq s\zeta(n, r), \qquad s := aD, \qquad r := a(S_0 + D + C) - 1.$$

This can be written in $(b, s_0)$ terms only, as

$$s = 2 + (2b - 1)/s_0,$$
$$r = \{2b(1 + s_0)^2 + 1 - 1/b\}/s_0.$$

When $b \geq 1$, the expression of $w_n^*$ stems from

$$\zeta(n, r) \leq \log(1 + n/r), \qquad r \geq 2b(1 + s_0)^2/s_0.$$

In the general case, $b \geq b_0 := 1/(1 + s_0)$ and $r(b, s_0) \geq r(b_0, s_0)$ since $r$ is increasing in $b$. Finally $r(b_0, s_0) = (2 + s_0)/s_0$ yields the value of $w_n^+$, since

$$\zeta(n, r) \leq \log(1 + n/r) \leq \log(1 + n/r(b_0, s_0)). \qquad \square$$

**6. General branching processes.** The results of Sections 3 and 4 can be extended, at a relatively low cost, to a wider context. Assume for instance that each particle $x$ in the $n$th generation gives birth to $Z_x \geq 1$ children, where $(Z_x)$ is i.i.d. and each $Z_x$ is distributed like $L_{n+1}$, say. On the event $\{Z_x = k\}$, order the $k$ children of $x$ from $y_1$ to $y_k$ and decide that the $k$-dimensional random vector $(\xi(y_1), \ldots, \xi(y_k))$ follows a given law $\pi_k^{n+1}$. Do this independently for different particles $x$ in the same generation and independently in different generations.

The PCR model is recovered when the law of $L_n$ is $(1 - \lambda_n)\delta_1 + \lambda_n\delta_2$ and when $\pi_1^n = \delta_0$ and $\pi_2^n = \delta_0 \otimes \pi$ for a given distribution $\pi$ on the nonnegative real numbers.

Coming back to the general setting, assume that every $L_n$ is integrable and call $z_n$ the size biased distribution of $L_n$, defined by $z_n(k) := k\alpha_k^n$, where

$$\alpha_k^n := \mathbb{P}(L_n = k)/\mathbb{E}(L_n).$$

Let the laws $\pi_k^n$ be (square) integrable for any $n \geq 1$ and any $k \geq 1$ in the support of the law of $L_n$. Let $\mu_k^n$ denote the expectation of $\xi(y_1) + \cdots + \xi(y_k)$ under $\pi_k^n$. Then the following analogue of Theorem A holds.

THEOREM J. *The law $\eta_n^*$ coincides with the law of the random variable*

$$\xi_1^* + \cdots + \xi_n^*,$$

*where $(\xi_n^*)_{n\geq 1}$ are independent and distributed as follows. For any fixed $n \geq 1$, draw $k \geq 1$ at random along the sized biased distribution $z_n$, then choose the index $i$ uniformly at random in $\{1, \ldots, k\}$ and let $\xi_n^*$ denote a copy of the $i$th marginal of $\pi_k^n$. Thus, for instance,*

$$\mathbb{M}(\eta_n^*) = \sum_{k=1}^{n} \sum_{j \geq 1} \mu_j^k \alpha_j^k.$$

In the PCR context, the only nonzero $\mu_j^k$ term is $\mu_2^k = \mu_k$ and $\alpha_2^k = \lambda_k/(1+\lambda_k)$ is $\alpha_k$. Thus $\mathbb{M}(\eta_n^*)$ is the sum of $\mu_k \alpha_k$ as in Theorem A.

The next step is to estimate the discrepancy between $\eta_n^*$ and $\eta_n$. With regard to first moments, their difference can now be negative or positive.

PROPOSITION 21. *One has*

$$\mathbb{M}(\eta_n) = \sum_{k=1}^{n} \sum_{j \geq 1} \mu_j^k \alpha_j^k (1 - \varepsilon_j^k),$$

*where the error terms $\varepsilon_j^k$, which can be positive or negative, are bounded by functions of $j$ and of the reproducing laws of $(L_i)_{i \leq k}$. Such bounds can be deduced from the inequalities*

$$0 \leq \mathbb{E}(L_k)\varepsilon_j^k - (j - \mathbb{E}(L_k))\mathbb{E}(1/S_{k-1})$$
$$\leq \mathbb{E}(\{(j - \mathbb{E}(L_k))^2 + (S_{k-1} - 1)\mathbb{V}(L_k)\}/S_{k-1}^2).$$

Assume, additionally, that the first moment of each marginal of $\pi_j^k$ is bounded by a given number $\mu_0^k$ (or that $|\mu_j^k| \leq j\mu_0^k$) for every $j \geq 1$ and $1 \leq k \leq n$. Then after some computations, Proposition 21 yields

$$|\mathbb{M}(\eta_n) - \mathbb{M}(\eta_n^*)| \leq \sum_{k=1}^{n} \mu_0^k \mathbb{E}(1/S_{k-1}) \mathbb{E}(L_k^3)/\mathbb{E}(L_k)^2.$$

Recalling finally that

$$\mathbb{E}(1/S_{k-1}) \leq \mathbb{E}(1/L_{k-1}) \cdots \mathbb{E}(1/L_1)/S_0,$$

we get the main result of this section.

THEOREM K. *Fix $n$ and fix some reproduction and mutation mechanisms for the $n$ first generations. Assume that there exists $r \geq 1$, $\mu_0$ and $L_0$ finite, such that, for every $1 \leq k \leq n$ and every $j \geq 1$,*

$$|\mu_j^k| \leq j^r \mu_0, \qquad \mathbb{E}(L_k^{2+r}) \leq L_0.$$

*Then there exists a finite constant $C := C(n, L_0)$ such that, for every $S_0 \geq 1$,*

$$|\mathbb{M}(\eta_n) - \mathbb{M}(\eta_n^*)| \leq C\mu_0 S_0^{-1}, \qquad \mathbb{M}(\eta_n^*) = \sum_{k=1}^{n} \sum_{j \geq 1} \mu_j^k \alpha_j^k.$$

*This holds with the (very crude) constant $C := nL_0$.*

We could deal with the second moment of $\eta_n$ along similar lines, but we leave this task to the interested reader.

**7. Estimation of mutation rates.** In the rest of the paper, $(\lambda_n)$ is a given deterministic sequence as in Section 5. Let $\widehat{\mu}$ denote the point estimator of $\mu$ by the first moment method, that is, the solution of the equation $t = \mathbb{E}(t)$ in the unknown $\mu$. Let $\widehat{\mu}_*$ denote its infinite-population limit, that is, the solution of $t = \mathbb{E}(t^*)$ in the unknown $\mu$. When the efficiency is constant, $\widehat{\mu}_*$ is the estimator due to Sun (1995).

COROLLARY 22. *We have $\widehat{\mu} > \widehat{\mu}_* = t/W_n$. More precisely, for any $S_0 \geq 1$,*

$$1 - \frac{r''}{S_0 + 1} \leq \frac{\widehat{\mu}_*}{\widehat{\mu}} \leq 1 - \frac{r}{S_0 + 1}, \qquad r := \frac{v_n}{W_n}, \qquad r'' := \frac{v_n''}{W_n}.$$

When is the finite-population correction to $\mathbb{E}(t^*)$ negligible? Assume that the $n$ first efficiencies $\lambda_k$ are greater than $\lambda^*$. Then (we omit the proof)

$$\left(1 - \frac{2}{n\lambda^* S_0}\right)\widehat{\mu} \leq \widehat{\mu}_* \leq \widehat{\mu}.$$

COROLLARY 23. *If $S_0 n \inf(\lambda_k) \gg 1$, then $\widehat{\mu}_* \approx \widehat{\mu}$.*

In the opposite direction, consider the situation where $n = 25$, $\lambda_k := \lambda_1/k$, $\lambda_1 := 0.25$ and $S_0 = 1$. Thus, the efficiency $\lambda_k$ decreases from $\lambda_1 = 0.25$ to $\lambda_{25} = 0.01$. Then $r = 0.495$ and $r'' = 0.770$, meaning that the true estimator $\widehat{\mu}$ is between $1.33 \cdot \widehat{\mu}_*$ and $1.63 \cdot \widehat{\mu}_*$.

With regard to second moment properties, we first make the following remark, direct from Theorem G.

COROLLARY 24. *In the infinite-population limit, $\widehat{\mu}_*$ is a consistent estimator of $\mu$. For finite populations, $\widehat{\mu}$ is not a consistent estimator of $\mu$.*

Assuming that $\widehat{\mu}_* \approx \widehat{\mu}$, the confidence interval of level $1 - 1/z^2$ for $\mu$ corresponds to $t$ lying in the interval bounded by

$$\mathbb{E}(t^*) \pm z\sqrt{\mathbb{V}(t^*)}.$$

In the Poisson case, $\mathbb{E}(\xi) = \mu = \mathbb{V}(\xi)$. Thus, replacing $\mu$ in $\mathbb{V}(t^*)$ by its point estimator $\widehat{\mu}_*$, we get an approximative confidence interval for $\mu$, bounded by the points $\widehat{\mu}_* \pm z\widehat{\sigma}_*/W_n$, where

$$\widehat{\mu}_* := t/W_n, \qquad \widehat{\sigma}_*^2 := (t + t^2 W_n'/W_n^2)/\ell.$$

When $\mu$ is small, $\widehat{\sigma}_*^2 \sim t/\ell$ and the interval is approximately bounded by

$$\widehat{\mu}_*(1 \pm z/\sqrt{t\ell}),$$

where $(t\ell)$ is the total number of mutations in the sample.

**8. Data set.** Weiss and von Haeseler (1997) applied their coalescent method to the data set of Saiki et al. (1988). The efficiency sequence $(\lambda_k)$ is not provided; neither is the initial population $S_0$. However, following Weiss and von Haeseler (1997), the extent of the amplification after 20, 25 and 30 cycles [provided by Saiki et al. (1988)] allows us to compute hypothetical efficiencies, which are constant and equal to $\lambda$ during the 20 first cycles, then constant and equal to $\lambda'$ during the cycles between 21 and 25, and finally constant and equal to $\lambda''$ during the cycles between 26 and 30. Numerically,

$$\lambda = 0.872, \qquad \lambda' = 0.743, \qquad \lambda'' = 0.146.$$

For a sample of size $\ell = 28$, 17 mutations were observed. Thus, $t = 17/28$. We find $W_{30} = 12.085$ and $\widehat{\mu}_* = t/W_{30} = 0.05024$. Furthermore, $v_{30}'' = 0.38435$ and $v_{30} = 0.03653$. Thus,

$$\widehat{\mu} \in \begin{cases} (0.05032, 0.05105), & \text{for } S_0 = 1, \\ (0.05025, 0.05039), & \text{for } S_0 = 10, \\ (0.05024, 0.05026), & \text{for } S_0 = 100. \end{cases}$$

These intervals are much smaller than the uncertainties associated to, first, the fact that the efficiencies are unknown and, second, the variation that a difference of $\pm 1$ in the count of the observed mutations would yield. With regard to the efficiencies, we followed Weiss and von Haeseler (1997) and chose a constant value from $k = 1$ to $k = 20$, then from $k = 21$ to $k = 25$, and finally from $k = 26$ to $k = 30$. More regular variations of $(\lambda_k)$ are possible.

The maximum likelihood method of Weiss and von Haeseler (1997) gives values that are similar to ours when $S_0$ is large and gives markedly different values when $S_0$ is small, namely an estimated value of $\mu$ of 0.060 for $S_0 = 1$. This suggests that the mean of the posterior distribution of $\mu$ is not the point where its density is maximal; in other words, that the distribution is far from being symmetric around its mean.

With regard to confidence intervals, $W_{30}' = 6.755$, and $\widehat{\sigma}_* = 0.149$, which is quite comparable to $\sqrt{t/\ell} = 0.147$; see the remark at the end of Section 7. This yields the interval bounded by $0.05024 \pm z\,0.01236$. For instance, $\mu \in (0.02552, 0.07496)$ at a level of confidence of 75%.

With regard to the variance, Weiss and von Haeseler (1997) simulate the correct value of 0.012 when $S_0$ is large and a variance of 0.016 when $S_0 = 1$.

## 9. On the mean of the sample.

9.1. *Boundary effects.* We recover some striking features of the numerical simulations in Weiss and von Haeseler (1995). Recall that the histogram of $\eta_n$ is always to the left of the histogram of $\eta_n^*$, meaning that $\eta_n^*$ stochastically dominates $\eta_n$ (this is specific to the case where $\xi \geq 0$ almost surely); see Piau (2002). Furthermore, the gap between the two distributions decreases to zero when $\lambda \to 1$. This last fact follows from $v_n' \leq 2(1-\lambda)$.

Another property that is not visible on the simulations is that, for $n$ fixed and when $\lambda \to 0$, the gap goes to zero as well. Our results prove that this effect appears only when $\lambda$ is so small that $n\lambda \ll 1$, that is, for values of the efficiency that Weiss and von Haeseler did not consider in their simulations.

9.2. *First generations effect.* Considering a smaller number of cycles in our test case in Section 7, we get similar values of the ratio $r$. For instance, if $n=5$, respectively, if $n=10$, then $r=0.521$, respectively, $r=0.516$. Roughly speaking, this means that the approximation fails only during the first generations, that is, when $S_k$ is not large enough yet.

9.3. *Estimating efficiencies.* For Bernoulli branching processes, the sequence $(\gamma_n S_n)$ is a positive martingale, bounded in $L^2$ (and in every $L^p$, $p \geq 2$). Thus, it converges, almost surely and in the mean to a random limit which is in $(0, +\infty)$ almost surely. This implies that

$$S_{n+1}/[S_n(1+\lambda_{n+1})] \to 1 \qquad \text{a.s.}$$

This fact, rather than the observation that $\mathbb{E}(S_{n+1}) = \mathbb{E}(S_n)(1+\lambda_{n+1})$, is the reason why $S_{n+1}/S_n$ is a good estimator of $(1+\lambda_{n+1})$. Additionally, some concentration of measure phenomena occurs; see the large deviations results of Athreya (1994) and Athreya and Vidyashankar (1995).

## 10. On the proof of stochastic recursions.

We first recall some basic tools. We have $\mathbb{M}(\zeta_n) = T_n/S_n$, where $T_n$ denotes the sum of the states of the population at time $n$. For any particle $x$ at time $n$, let $\varepsilon(x) := 1$ if $x$ has two children; otherwise, let $\varepsilon(x) := 0$. Let $(\xi(x))$ denote i.i.d. random variables distributed like $\xi$. Then

$$T_{n+1} = \sum_{x=1}^{S_n} s(x)(1+\varepsilon(x)) + \xi(x)\varepsilon(x), \qquad S_{n+1} = \sum_{x=1}^{S_n} 1 + \varepsilon(x).$$

The exchangeability of $\varepsilon(x)$ implies that

$$S_n \mathbb{E}\left(\frac{1+\varepsilon(x)}{S_{n+1}} \Big| \mathcal{F}_n\right) = 1.$$

Likewise,
$$S_n \mathbb{E}\left(\frac{\varepsilon(x)}{S_{n+1}}\Big|\mathcal{F}_n\right) = 1 - H(S_n).$$

Integrating $T_{n+1}/S_{n+1}$, we get
$$\mathbb{E}(\mathbb{M}(\zeta_{n+1})|\mathcal{F}_n) = \mathbb{M}(\zeta_n) + \mu[1 - H(S_n)].$$

Turning to the evaluation of the second moment, we use once again the expression of $T_{n+1}$ and the exchangeability properties that arise. Separating carefully the square terms from the rectangular terms in $T_{n+1}^2$, and skipping the details, we get

$$\mathbb{E}(\mathbb{M}(\zeta_{n+1})^2|\mathcal{F}_n) = \mathbb{M}_2(\zeta_n)B''(S_n) + \mathbb{M}(\zeta_n)^2 B_1(S_n)$$
$$+ 2\mu \mathbb{M}(\zeta_n)[1 - H(S_n)] + \nu B(S_n) + \mu^2 B_2(S_n),$$

where the only coefficients that are not already defined are
$$B_1(1) := 0, \qquad B_1(k) := k^2 \mathbb{E}(L_1 L_2/M_k^2), \qquad k \geq 2,$$
$$B_2(k) := \mathbb{E}((1 - k/M_k)^2).$$

We take the expectation of both sides and substract to this the recursion relationship for $\mathbb{E}(\mathbb{M}(\zeta_n))$ squared. The $\mu$ terms cancel. The $\mu^2$ terms add to
$$B_2(S_n) - [1 - H(S_n)]^2 = B'(S_n).$$

The $\nu$ term is $B(S_n)$. The $\mathbb{M}_2(\zeta_n)$ and $\mathbb{M}(\zeta_n)^2$ terms remain and we must show that the coefficients of these, namely $B''(S_n)$ and $B_1(S_n) - 1$, add to 0. Since $M_k^2$ is the sum of $k$ squares $L_i^2$ and of $k(k-1)$ products $L_i L_j$ for $i \neq j$, the exchangeability of $(L_i)$ yields
$$k\mathbb{E}((L_1^2 - L_1 L_2)/M_k^2) + k^2 \mathbb{E}(L_1 L_2/M_k^2) = 1,$$

that is, $B''(k) + B_1(k) = 1$. This ends the proof of Theorem E. The proof of Lemma 11 uses similar techniques and we omit it.

**11. Additional term—proofs.** We begin with the proof of Lemma 9 and with simple considerations about the sequences $B$, $B'$ and $B''$ that are valid for any distribution of $L$. Laplace's method and the law of large numbers yield that, when $k \to \infty$,

$$kB(k) \to \frac{\mathbb{E}(L) - 1}{\mathbb{E}(L)^2}, \qquad kB'(k) \to \frac{\mathbb{V}(L)}{\mathbb{E}(L)^4}, \qquad kB''(k) \to \frac{\mathbb{V}(L)}{\mathbb{E}(L)^2}.$$

This implies that Lemma 9 cannot hold with $b < \alpha(1-\alpha)$, $b' < \lambda(1-\lambda)/(1+\lambda)^4$ or $b'' < \lambda(1-\lambda)/(1+\lambda)^2$. Thus, Lemma 9 is optimal when $\lambda \to 0$.

PROOF OF LEMMA 9. The first assertion of Lemma 9 stems from the concavity of $(m-k)/m^2$ with respect to $m$ on the interval $(k, 2k)$. The $B''$ result follows from the facts that, for $k \geq 2$,

$$B''(k) = \lambda(1-\lambda)\mathbb{E}\left(\frac{k}{(3+M_{k-2})^2}\right),$$

that $M_{k-2} \geq k-2$ a.s. and that $k/(k+1)^2 \leq 1/(k+2)$. We postpone (a stronger version of) the $B'$ assertion to Lemma 25. □

REMARK. In the Bernoulli case, $B''(k) \geq (\mathbb{V}(L)/\mathbb{E}(L)^2) \cdot k/(k+1)^2$.

REMARK. We can refine the uniform bounds as follows. First, $v_n \leq 1$ and $v_n'' \leq 3$. Assuming that $\lambda_k \geq \lambda^*$ for every $k \leq n$,

$$v_n \leq 1 - \lambda^*.$$

On the other hand, assuming that $\lambda_k \leq \lambda^+$ for every $k \leq n$,

$$v_n \geq (1-\gamma_n)(1-\lambda^+).$$

This yields lower bounds of $Z_n$ and $R_n$ since, for instance,

$$\nu v_n/(S_0+1) \leq Z_n \leq (\nu+\mu^2)v_n/(S_0-1).$$

Lemma 25 is a key step in the proof of Lemma 26.

LEMMA 25. (i) *One has* $G(k) \leq 1/(1+\lambda)^2 + 3A(k)$.
(ii) *For any integer* $j \geq 1$,

$$\mathbb{E}((k/M_k)^j) \leq 1/(1+\lambda)^j + A(k)j(j+1)/2.$$

(iii) *One has* $B'(k) \leq A(k)(1+3\lambda)/(1+\lambda)$.
(iv) *One has* $B'(k) \leq b'/(k+1)$, *with*

$$b' := \lambda(1-\lambda)(1+3\lambda)/(1+\lambda)^2 \leq \lambda.$$

PROOF. Assertion (iv) follows from assertion (iii) and from the upper bound of $A(k)$ in Corollary 34. Assertion (iii) follows from assertion (i) and from the fact that

$$B'(k) = G(k) - H(k)^2 \leq [3 - 2/(1+\lambda)]A(k).$$

The convexity of the function $m \mapsto 3k/m - k^2/m^2$ on $m \geq k$ implies that

$$3\mathbb{E}(k/M_k) - \mathbb{E}((k/M_k)^2) \geq 3/(1+\lambda) - 1/(1+\lambda)^2.$$

Going back to the definitions, this is assertion (i). For assertion (ii), we use the convexity of the function $m \mapsto (j+2)/m^j - jk/m^{j+1}$ on $m \geq k$ to perform a recursion on $j$. □

The contribution of $\nu$, respectively, of $\mu^2$, to $R_{n+1} - R_n$ is greater than $\mathbb{E}(B(S_n))$, respectively, $\mathbb{E}(B'(S_n))$. The contribution of $\nu$, respectively, of $\mu^2$, to $Z_{n+1} - Z_n$ is $\mathbb{E}(A(S_n))$, respectively, $A'_{n+1}$. From Corollary 34, $A(k) \leq \alpha(1-\lambda)/2$. Hence,

$$A'_{n+1} \leq \mathbb{E}(A'(S_n)), \qquad \text{where } A'(k) := (1-\lambda)A(k).$$

Thus, Lemma 26 below implies that $R_n \geq Z_n/2$, that is, the $\ell \geq 3$ point in Proposition 19.

LEMMA 26. *For any $k \geq 1$, $B(k) \geq A(k)/2$ and $B'(k) \geq A'(k)/2$.*

PROOF. The definitions and two lines of algebra show that the inequality $B(k) \geq A(k)/2$ is equivalent to

$$(k-2)A(k) + 2G(k) \leq 2(1-\alpha).$$

From assertion (i) of Lemma 25, a sufficient condition is

$$(k+4)A(k) \leq 2\alpha(1-\alpha).$$

The upper bound of $A(k)$ in Corollary 34 implies that this holds if

$$(k+4)/(k+1) \leq 2/(1-\lambda^2).$$

This settles the $k \geq 2$ case, since $(k+4)/(k+1) \leq 2$ for $k \geq 2$. Because the $k = 1$ case is obvious, this proves that $B(k) \geq A(k)/2$ for any $k \geq 1$.

Our proof of $B'(k) \geq A'(k)/2$ is more intricate. According to Corollary 36 in Section 12.4, whose notation we use from now on, it is enough to check that

$$\widetilde{G}(k) \geq \widetilde{H}(k)^2 + (1-\lambda)(1+\lambda)(\widetilde{H}(k)-1)/2,$$

that is, after some rearrangements, to check that

$$(1+\lambda^2)G_1 k \geq (3+\lambda^2)G_2 + (1-\lambda^2)\frac{G_3}{k} + 2\left(G_1^2 + \frac{G_2^2}{k^2} + \frac{G_3^2}{k^4}\right).$$

First, $G_2 \leq (1-2\lambda)G_1$. Furthermore, since $n_2 \in (0, 1/4)$, for any $k \geq 2$,

$$G_3 = G_1[1 + 3(k-2)n_2]/(1+\lambda) \leq G_1 uk, \qquad u := 3/4.$$

For $k \geq 2$, this yields the sufficient condition

$$(1+\lambda^2)k \geq (3+\lambda^2)(1-2\lambda) + u(1-\lambda^2)$$
$$+ 2G_1 + 2\frac{G_1}{k^2}((1-2\lambda)^2 + u^2).$$

Since $G_1 \leq \lambda(1-\lambda)$, the sum of the three first terms on the right-hand side is the sum of $3 + u$ and of a polynomial in $\lambda$ with negative coefficients, hence at most $3 + u$. For any $k \geq 4$, the last term on the right-hand side is at most $2G_1(1+u^2)/k^2 \leq 1/16$ because $u \leq 1$, $G_1 \leq 1/4$ and $k \geq 4$. Hence, the right-hand side is at most $3 + 3/4 + 1/16 < 4$. Since the left-hand side is at least $k$, we are done for any $k \geq 4$.

Finally, setting $D(k, \lambda) := (2B'(k) - A'(k))(1+\lambda)/[\lambda(1-\lambda)]$, we get, with the help of Maple$^{©}$ software,

$$D(1, \lambda) = \lambda \geq 0,$$
$$18D(2, \lambda) = 2 + 10\lambda - 9\lambda^2 + \lambda^3 \geq 2,$$
$$200D(3, \lambda) = 24 + (1-\lambda)(1 + 58\lambda - 55\lambda^2 - \lambda^4) \geq 24.$$

Hence, $B'(k) \geq A'(k)/2$ for any $k \geq 1$. □

## 12. Harmonic moments—statements.

12.1. *Method.* Following the technique of Piau (2002, 2004a), we seek to compare $A(k)$ with $1/k$ and to bound $H_y(k)$. Iterating these bounds yields good estimations of the harmonic moments of $S_k$.

A remarkable feature of the rescaled harmonic mean $H_y(k)$ of the sum of $k$ i.i.d. positive random variables is that, for $y \geq 0$, $H_y(k)$ is a decreasing function of $k$ for $k \geq 1$. This very general fact seems to be unknown. It holds in a wider context (see Lemma 30) and we prove it by a completely deterministic method in Section 13.1.

In the restricted context of the Bernoulli branching process, we uncover two other monotonicities. First, for $y \geq 1$, $H_{-y}(k)$ is an increasing function of $k$ for $k > y$. Second, a suitably normalized correction of $H(k)$ is decreasing (however, see Remark in Section 12.3). Thus, Lemmas 30, 31 and 33 are crucial steps in our proofs. We state the following consequence.

COROLLARY 27. *In the Bernoulli case, for any $y \geq 0$ and $S_0 \geq 1$,*

$$\gamma_n/(S_0 + y) \leq \mathbb{E}[1/(S_n + y)] \leq \gamma_n^{(y+2)}/(S_0 + y).$$

*For any $S_0 > y \geq 1$,*

$$\mathbb{E}[1/(S_n - y)] \leq \gamma_n/(S_0 - y).$$

*As a consequence, for any $S_0 \geq 2$,*

$$\gamma_n/S_0 \leq \mathbb{E}(1/S_n) \leq \gamma_n/(S_0 - 1).$$

REMARK 28. If $S_0 \geq 2$, $\mathbb{E}(1/S_n)$ is exactly of order $\gamma_n$, that is, of the order of $1/\mathbb{E}(S_n)$. If $S_0 = 1$, our results show only that $\mathbb{E}(1/S_n)$ is at most of order $\gamma_n^{(2)}$, which can be much greater than $\gamma_n = 1/\mathbb{E}(S_n)$. An upper bound of $\mathbb{E}(1/S_n)$ of order $\gamma_n$ indeed holds when $S_0 = 1$, namely

(5) $$\mathbb{E}(1/S_n) \leq \gamma_n(1 + 1/\lambda_n^*),$$

where $\lambda_n^*$ is the minimum of the sequence $(\lambda_k)$ up to time $n$ (we omit the proof). This bound is optimal, up to a factor 2, since, when $\lambda_k = \lambda$ is constant, we can show that the limit of $\mathbb{E}(1/S_n)/\gamma_n$ is at least $(1 + 1/\lambda)/2$; see Piau (2004a or b).

Equation (5) could be used to replace the products $\gamma_n^{(2)}$ in all our upper bounds by $\gamma_n(1 + 1/\lambda_n^*)$. Recall finally that, in actual reactions, the sequence $(\lambda_k)$ is nonincreasing. Then $\lambda_n^* = \lambda_n$ and $\gamma_n(1 + 1/\lambda_n^*) = \gamma_{n-1}/\lambda_n$.

12.2. *Monotonicities of $H_y^\varphi$.*

DEFINITION 29. For any integer $k \geq 1$ and real number $y$ such that $k + y > 0$, and for any nonnegative function $\varphi$, define
$$H_y^\varphi(k) := \mathbb{E}\left(\varphi\left(\frac{M_k + y}{k + y}\right)\right).$$

We recover $H_y$ when $\varphi$ is the convex and decreasing function $\varphi(x) = 1/x$. For any locally bounded $\varphi$ and any $L \geq 1$ of bounded support, the law of large numbers and an easy domination imply that $H_y^\varphi(k) \to \varphi(\mathbb{E}[L])$ when $k \to \infty$.

LEMMA 30. (i) *Assume that $(L_i)$ is exchangeable and that $\varphi$ is convex. Then $H_0^\varphi$ is nonincreasing. If furthermore $(L_i)$ is i.i.d., then $H_0^\varphi(k) \geq H_0^\varphi(\infty) = \varphi(\mathbb{E}[L])$.*

(ii) *Assume that $y \geq 0$, that $(L_i)$ is exchangeable, and that $\varphi$ is convex and nonincreasing. Then $H_y^\varphi$ is nonincreasing. If $(L_i)$ is furthermore i.i.d., then $H_y^\varphi(k) \geq H_y^\varphi(\infty) = \varphi(\mathbb{E}[L])$.*

For instance, $H_y$ decreases from $H_y(1) = \mathbb{E}[(1+y)/(L+y)]$ to $1/\mathbb{E}(L)$ for any $y \geq 0$. For Bernoulli branching processes, we get
$$1 - \alpha \leq H_y(k) \leq 1 - \lambda/(y+2).$$
Furthermore, $H_y$ describes one step of the evolution of $1/(S_n + y)$, since
$$H_y(S_n) = \mathbb{E}\left(\frac{S_n + y}{S_{n+1} + y}\Big|S_n\right) \qquad \text{a.s.}$$
This yields the first part of Corollary 27.

LEMMA 31. *For Bernoulli branching processes and $y \geq 1$, the sequence $H_{-y}(k)$ is increasing for $k > y$. Thus,*
$$H_{-y}(k) \leq 1/\mathbb{E}(L) = 1 - \alpha.$$

Iterating Lemma 31 yields the second part of Corollary 27.

REMARK. For non-Bernoulli laws, $H_{-1}$ can fail to be increasing.

12.3. *Monotonicity of A.* Lemmas 32 and 33 are proved in Section 13.2.

LEMMA 32. *For $k=1$, $A(1) = \mathbb{E}(1/L) - 1/\mathbb{E}(L)$. When $k \to \infty$,*
$$kA(k) \to \mathbb{V}(L)/\mathbb{E}(L)^3.$$

LEMMA 33. *For Bernoulli branching processes, $(k+1)A(k)$ is decreasing.*

For Bernoulli branching processes, $A(1) = \alpha(1-\lambda)/2$ and $\mathbb{V}(L) = \lambda(1-\lambda)$.

COROLLARY 34. *For Bernoulli branching processes,*
$$\alpha(1-\lambda)/(1+\lambda)^2 \leq (k+1)A(k) \leq \alpha(1-\lambda).$$

COROLLARY 35. *In particular, $kA(k) \leq \alpha(1-\lambda)$.*

REMARK. (i) The sequence $kA(k)$ is not always decreasing, even in the restricted case of Bernoulli branching processes. If $\lambda \leq 2 - \sqrt{3}$, we can show that $kA(k)$ is in fact increasing (we omit the proof).

(ii) In the general context, the sequence $(k+1)A(k)$ is not always decreasing. If the law of $L$ is $(1-p)\delta_1 + p\delta_3$ and if $p \leq 1/16$, we can show that $(k+1)A(k)$ is in fact increasing (we omit the proof).

12.4. *Taylor expansions.* Precise estimates of the remaining terms of the Taylor expansions of the functions $x \mapsto 1/x$ and $x \mapsto 1/x^2$ yield bounds of $G(k)$ and $H(k)$, hence of $A(k)$. Set $n_2 := \lambda(1-\lambda)$. The following results are proved in Section 13.3.

COROLLARY 36. *One has $G(k) \geq \widetilde{G}(k)/(1+\lambda)^2$ and $H(k) \leq \widetilde{H}(k)/(1+\lambda)$ with*
$$\widetilde{G}(k) := 1 + 3G_1/k - 4G_2/k^2 + 2G_3/k^3,$$
$$\widetilde{H}(k) := 1 + G_1/k - G_2/k^2 + G_3/k^3,$$
*where*
$$G_1 := \frac{n_2}{(1+\lambda)^2}, \qquad G_2 := \frac{n_2(1-2\lambda)}{(1+\lambda)^3}, \qquad G_3 := \frac{n_2(1+3(k-2)n_2)}{(1+\lambda)^3}.$$

Corollaries 34 and 36 provide tight bounds of $A(k)$. For instance, when $\lambda \to 0$, $A(k)/\lambda$ is asymptotically between $1/k - 1/k^2 + 1/k^3$ and $1/(k+1)$. The difference is $1/[k^3(k+1)]$. In the general case, the following corollary holds.

COROLLARY 37. *For any $k \geq 1$,*
$$\frac{\lambda(1-\lambda)}{(1+\lambda)^3}\frac{1}{k+1} \leq A(k) \leq \frac{\lambda(1-\lambda)}{(1+\lambda)^3}\frac{k+1}{k^2}.$$

*Therefore, for any $k \geq 2$,*
$$A(k) \leq \frac{\lambda(1-\lambda)}{(1+\lambda)^3}\frac{1}{k-1}.$$

## 13. Harmonic moments—proofs.

13.1. *Proof of the monotonicities of $H_y^\varphi$.* Surprisingly (to us), Lemmas 30 and 31 reflect almost sure properties that have nothing to do with randomness. Assume first that $y = 0$ and note that $M_{k+1}/(k+1)$ is the barycenter with equal weights of the $(k+1)$ random variables $M_k^{(i)}/k$, where $M_k^{(i)} := M_{k+1} - L_i$. Thus, the convexity of $\varphi$ implies that

(6) $$\varphi\left(\frac{M_{k+1}}{k+1}\right) \leq \frac{1}{k+1}\sum_{i=1}^{k+1}\varphi\left(\frac{M_k^{(i)}}{k}\right).$$

By exchangeability, each $M_k^{(i)}$ is distributed like $M_k$. Taking the expectations of both sides of (6) yields the result.

For $y \geq 0$, set $N_k := (M_k + y)/(k+y)$ and $y_k := y/[k(y+k+1)] \geq 0$. Tedious computations show that

$$N_{k+1} = -y_k + (1+y_k)\frac{1}{k+1}\sum_{i=1}^{k+1}N_k^{(i)},$$

where $M_k^{(i)}$ yields $N_k^{(i)}$ like $M_k$ yields $N_k$. Since each $N_k^{(i)} \geq 1$, $N_{k+1}$ is greater than the barycenter of the random variables $N_k^{(i)}$, which are distributed like $N_k$. Since $\varphi$ is nonincreasing, taking expectations yields Lemma 30.

Likewise, the proof of Lemma 31 for Bernoulli random variables is entirely nonrandom. Assume that $L_j = 2$ for $i$ indices $j$ between 1 and $k+1$, and assume that $L_j = 1$ for the $k+1-i$ other indices $j$. Then $M_{k+1} = k+i+1$, $M_k^{(j)} = k+i-1$ for $i$ indices $j$ and $M_k^{(j)} = k+i$ for the $k+1-i$ other indices $j$. Thus, it is enough to check the almost sure inequality

$$\frac{k+1-y}{k+i+1-y} \geq \frac{k-y}{k+1}\left(\frac{i}{k+i-1-y} + \frac{k+1-i}{k+i-y}\right)$$

for any $0 \leq i \leq k+1$ and $k > y$. After simplifications, this is equivalent to the condition that either $i = 0$ or $i \geq 1$ and

$$(y-1)k + (y+1)i - 1 - y^2 \geq 0,$$

which holds for any $i \geq 1$ and $k > y \geq 1$. Thus, Lemma 31 holds.

13.2. *Proof of the properties of A.* For $|t| \leq 1$, let $f(t) := \mathbb{E}(t^L)$ and $g(t) := \mathbb{E}((L_1 - L_2)^2 t^{L_1 + L_2 - 3})/2$. An integration by parts yields

$$A(k) = \int_0^1 g(t) f'(t)^{-2} f(t)^k \, dt.$$

(We omit the details.) Since $|f(t)| < 1$ for $|t| < 1$ and $f(1) = 1$, when $k \to \infty$, the integral is controlled by the behavior of the integrand when $t \to 1$. More precisely, from Laplace's method,

$$A(k) \sim g(1) f'(1)^{-2} / [k f'(1)],$$

where $g(1)$ and $f'(1)$ denote limits when $t \to 1$, $t < 1$. Since $g(1) = \mathbb{V}(L)$ and $f'(t) \to \mathbb{E}(L)$, Lemma 32 holds.

For Bernoulli branching processes, $g(t) = \lambda(1 - \lambda)$ is constant. Thus, $A(k)$ takes the simpler form $A(k) = \lambda(1 - \lambda) I(k, 1)$, where

$$I(k, \ell) := \int_0^1 (f)^k (f')^{-2\ell}.$$

The facts that $(k+1) f' f^k$ is the derivative of $f^{k+1}$ and that $f'' = 2\lambda$, together with an integration by parts yield the recursion

$$(k+1) I(k, \ell) = (1 + \lambda)^{-(2\ell+1)} + 2\lambda(2\ell + 1) I(k+1, \ell+1).$$

Since $I(\cdot, \ell)$ is decreasing for any $\ell$, this implies that the sequence $(k+1) I(k, \ell)$ is a decreasing function of $k$ and, finally, that Lemma 33 holds.

13.3. *Proofs of Taylor expansions.* Let $T_i(f)(x_0, \cdot)$ denote the Taylor expansion at $x_0$ of the function $f$, up to order $i \geq 0$. For instance, if $h(x) := 1/x$ and $g(x) := 1/x^2$,

$$T_i(h)(x_0, x) = \sum_{j=0}^i (-1)^j (x - x_0)^j / x_0^{j+1},$$

$$T_i(g)(x_0, x) = \sum_{j=0}^i (-1)^j (j+1)(x - x_0)^j / x_0^{j+2}.$$

We now estimate the remaining terms, perhaps more precisely than is usual.

LEMMA 38. *For any $i \geq 1$,*

$$h(x) = T_{2i-1}(h)(x_0, x) + (x - x_0)^{2i} / (x_0^{2i} x),$$
$$g(x) = T_{2i-1}(g)(x_0, x) + (x - x_0)^{2i} (x_0 + 2ix) / (x_0^{2i+1} x^2).$$

PROOF. Recursion on $i \geq 1$. □

COROLLARY 39. *If $x$ and $x_0$ are both in a positive interval $(x_-, x_+)$,*

$$h(x) \leq T_{2i-1}(h)(x_0, x) + h_i(x - x_0)^{2i},$$
$$g(x) \geq T_{2i-1}(g)(x_0, x) + g_i(x - x_0)^{2i},$$

*where*

$$h_i := \frac{1}{x_0^{2i} x_-}, \qquad g_i := \frac{2ix_+ + x_0}{x_0^{2i+1} x_+^2} \geq \frac{2i}{x_0^{2i+1} x_+}.$$

*A similar lower bound of $h(x)$ and a similar upper bound of $g(x)$ hold.*

REMARK. Taylor–Lagrange formulas yield greater error terms for $h(x)$ and for $g(x)$ that are, respectively,

$$h_i' := \frac{1}{x_-^{2i+1}} \quad \text{and} \quad g_i' := \frac{2i+1}{x_+^{2i+2}}.$$

We apply Corollary 39 to $x := M_k/k$, $x_0 := 1 + \lambda$, $x_- := 1$ and $x_+ := 2$, and $i := 2$. Introducing $m_j := \mathbb{E}((x - x_0)^j)$ and using $m_1 = 0$, we get

$$H(k) \leq 1/(1+\lambda) + m_2/(1+\lambda)^3 - m_3/(1+\lambda)^4 + m_4/(1+\lambda)^4,$$
$$G(k) \geq 1/(1+\lambda)^2 + 3m_2/(1+\lambda)^4 - 4m_3/(1+\lambda)^5 + 2m_4/(1+\lambda)^5.$$

If $n_j$ denotes the $j$th moment of $L - \mathbb{E}(L)$, $m_2 = n_2/k$, $m_3 = n_3/k^2$ and

$$m_4 = (n_4 + 3(k-1)n_2^2)/k^3.$$

If $L$ is Bernoulli, $n_2 = \lambda(1-\lambda)$, $n_3 = n_2(1-2\lambda)$ and $n_4 = n_2(1-3n_2)$. This yields Corollary 36.

PROOF OF COROLLARY 37. The expression of $\widetilde{H}(k)$ in Corollary 36 yields

$$kA(k) \leq \frac{\lambda(1-\lambda)}{(1+\lambda)^3}\left(1 + \frac{a}{k}\right),$$

where, after some rearrangements, $a$ can be written as

$$a = \lambda - \left(1 - \frac{1}{k}\right)\frac{1 - 4\lambda(1-\lambda)}{1+\lambda} - \frac{2\lambda(1-\lambda)}{(1+\lambda)k}.$$

Thus, $a \leq 1$ and the upper bound is proved. The lower bound is in Corollary 34. □

# REFERENCES

Athreya, K. B. (1994). Large deviation rates for branching processes. I. Single type case. *Ann. Appl. Probab.* **4** 779–790. MR1284985

Athreya, K. B. and Vidyashankar, A. N. (1995). Large deviation rates for branching processes. II. The multitype case. *Ann. Appl. Probab.* **5** 566–576. MR1336883

Higuchi, R., Dollinger, G., Walsh, P. S. and Griffith, R. (1992). Simultaneous amplification and detection of specific DNA sequences. *Biotechnology* **10** 413–417.

Higuchi, R., Fockler, C., Dollinger, G. and Watson, R. (1993). Kinetic PCR analysis: Real-time monitoring of DNA amplification reactions. *Biotechnology* **11** 1026–1030.

Jagers, P. and Klebaner, F. (2003). Random variations and concentration effects in PCR. *J. Theoret. Biol.* **224** 299–304. MR2067238

Piau, D. (2002). Mutation-replication statistics for polymerase chain reactions. *J. Comput. Biol.* **9** 831–847.

Piau, D. (2004a). Immortal branching Markov processes: Averaging properties and polymerase chain reaction applications. *Ann. Probab.* **32** 337–367. MR2040785

Piau, D. (2004b). Harmonic moments of branching processes. Preprint.

Saiki, R. K., Gelfand, D. H., Stoffel, S., Scharf, S. J., Higuchi, R., Horn, G. T., Mullis, K. B. and Erlich, H. A. (1988). Primer-directed enzymatic amplification of DNA with a thermostable DNA polymerase. *Science* **239** 487–491.

Schnell, S. and Mendoza, C. (1997). Enzymological considerations for a theoretical description of the quantitative competitive polymerase chain reaction (QC-PCR). *J. Theoret. Biol.* **184** 433–440.

Sun, F. (1995). The polymerase chain reaction and branching processes. *J. Comput. Biol.* **2** 63–86.

Wang, D., Zhao, C., Cheng, R. and Sun, F. (2000). Estimation of the mutation rate during error-prone polymerase chain reaction. *J. Comput. Biol.* **7** 143–158.

Weiss, G. and von Haeseler, A. (1995). Modeling the polymerase chain reaction. *J. Comput. Biol.* **2** 49–61.

Weiss, G. and von Haeseler, A. (1997). A coalescent approach to the polymerase chain reaction. *Nucleic Acids Research* **25** 3082–3087.

ex-LaPCS—Domaine de Gerland
Université Claude Bernard Lyon 1
50 avenue Tony-Garnier
69366 Lyon Cedex 07
France
e-mail: didier.piau@univ-lyon1.fr
url: http://lapcs.univ-lyon1.fr